\documentclass[11pt]{article}
\usepackage{amsmath}
\usepackage{mathrsfs}
\usepackage{amsfonts}

\usepackage{amsfonts, amsmath, amssymb}
\usepackage{amssymb,amsfonts,amsmath,
latexsym, epsfig,cite, psfrag,eepic,colordvi}
\usepackage{amscd,graphics}
\usepackage{lineno}

\newcommand{\qed}{\hfill $\Box $}
\newcommand{\pf}{\noindent {\bf Proof.} }

\newtheorem{theorem}{Theorem}[section]
\newtheorem{lemma}[theorem]{Lemma}
\newtheorem{coro}[theorem]{Corollary}

\begin{document}

\title{Maximum spectral radius of graphs
 with given connectivity and minimum degree}

\author{Hongliang Lu\,\textsuperscript{1}, \
Yuqing Lin\,\textsuperscript{2}
\\ {\small $^1$ Department of Mathematics}
\\ {\small Xi'an Jiaotong University, Xi'an 710049, PR China}
\\ {\small $^2$ School of Electrical Engineering and Computer Science}
\\ {\small The University of Newcastle, Newcastle, Australia}
}

\date{}

\maketitle

\date{}

\maketitle

\begin{abstract}
Shiu, Chan and Chang [On the spectral radius of graphs with
connectivity at most $k$, J. Math. Chem., 46 (2009), 340-346]
studied the spectral radius of graphs of order $n$ with $\kappa(G)
\leq k$ and showed that among those graphs, the maximum spectral
radius is obtained uniquely at $K_k^n$, which is the graph obtained
by joining $k$ edges from $k$ vertices of $K_{n-1}$ to an isolated
vertex. In this paper, we study the spectral radius of graphs of
order $n$ with $\kappa(G)\leq k$ and minimum degree $\delta(G)\geq k
$. We show that among those graphs, the maximum spectral radius is
obtained uniquely at $K_{k}+(K_{\delta-k+1}\cup K_{n-\delta-1})$.

\end{abstract}

\begin{flushleft}
{\em Key words:} connectivity; spectral radius
\\
\end{flushleft}

\section{Introduction}

Let $G$ be a simple graph of order $n$ with vertex set $V(G) =
\{v_1, v_2, \ldots, v_n\}$. We denote by $\delta (G)$ the minimum
degree of   vertices of  $G$ . The adjacency matrix of the graph $G$
is defined to be a matrix $A(G) = [a_{ij}]$ of order $n$, where
$a_{ij} = 1$ if $v_i$ is adjacent to $v_j$, and $a_{ij} = 0$
otherwise. Since $A(G)$ is symmetric and real, the eigenvalues of
$A(G)$, also referred to as the eigenvalues of $G$, can be arranged
as: $\lambda_n(G)\leq \lambda_{n-1}(G)\leq\cdots \leq \lambda_1(G)$.
The largest eigenvalue $\lambda_1(G)$ is called spectral radius and
also denoted by $\rho(G)$. For $k\geq 1$, we say that a graph $G$ is
$k$-connected if either $G$ is a complete graph $K_{k+1}$, or else
it has at least $k+2$ vertices and contains no $(k-1)$-vertex cut.
The \emph{connectivity} $\kappa(G)$ of $G$ is the maximum value of
$k$ for which $G$ is $k$-connected.

When $G$ is connected, $A(G)$ is irreducible and by the
Perron-Frobenius Theorem, the spectral radius is simple and there is
an unique positive unit eigenvector. We shall refer to such an
eigenvector as the Perron vector of $G$.

 The eigenvalues of a graph are related to many of its
properties and key parameters. The most studied eigenvalues have
been the spectral radius $\rho(G)$ (in connection with the chromatic
number, the independence number and the clique number of the graph
\cite{Nikiforov,Wilf}). 
Brualdi and Solheid \cite{Brualdi} proposed the following problem
concerning spectral radius:

\emph{Given a set of graphs $\mathscr{C}$, find an upper bound for
the spectral radius of graphs in $\mathscr{C}$ and characterize the
graphs in which the maximal spectral radius is attained.}

If $\mathscr{C}$ is the set of all connected graphs on $n$ vertices
with $k$ cut vertices, Berman and Zhang \cite{Berman} solved this
problem. Liu et al. \cite{Liu} studied this problem for
$\mathscr{C}$ to be the set of all  graphs on $n$ vertices with $k$
cut edges. Wu et al. \cite{Wu} studied this problem for
$\mathscr{C}$ to be the set of trees on $n$ vertices with $k$
pendent vertices. Feng, Yu and Zhang \cite{Feng} studied this
problem for $\mathscr{C}$ to be the set of graphs on $n$ vertices
with matching number $\beta$.

Li, Shiu, Chan and Chang   studied this question for graphs with $n$
vertices and connectivity at most $k$, and obtained the following
result.

\begin{theorem}[Li, Shiu, Chan and Chang \cite{Shiu}]
Among all the graphs with connectivity at most $k$, the maximum
spectral radius is obtained uniquely at $K_k+(K_1\cup K_{n-k-1})$.
\end{theorem}

Let $G_{k,\delta,n}=K_k+(K_{\delta-k+1}\cup K_{n-\delta-1})$. We
denote by $\mathcal {V}_{k,\delta,n}$ the set of graphs of order $n$
with $\kappa(G) \leq k \leq n-1$ and $\delta(G)\geq k$. Clearly,
$\mathcal {V}_{k,\delta+1,n}\subseteq \mathcal {V}_{k,\delta,n}$. In
this paper, we investigate the problem for the graphs in $\mathcal
{V}_{k,\delta,n}$. We show that among all those graphs, the maximal
spectral radius is obtained uniquely at $G_{k,\delta,n}$.

In our arguments, we need the following   technical lemma.
\begin{theorem}[Li and Feng  \cite{Li}]\label{lem2}
Let $G$ be a connected graph, and $G'$ be a proper subgraph of $G$.
Then $\rho(G') < \rho(G)$.
\end{theorem}

\section{Main results}

\begin{theorem}\label{thm1}
Let $u$, $v$ be two vertices of the connected graph $G$. Let
$\{v_1,\ldots, v_k\}\subseteq N(v)$   and
$\{v_{k+1},\ldots,u_{k+l}\}\subseteq V(G)-N(v)$. Suppose $x = (x_1,
x_2,\ldots, x_n)^T$ is the Perron vector of $G$, where $x_i$
corresponds to the vertex $v_i (1\leq i\leq  n)$. Let $G^*$ be the
graph obtained from $G$ by deleting the edges $vv_i\ (1\leq i\leq
k)$ and adding the edges $vv_i\ (k+1\leq i\leq k+l)$. If
$\sum_{i=1}^{k}x_i\leq \sum_{i=k+1}^{k+l}x_i$, then $\rho(G)\leq
\rho(G^*)$. Furthermore, if $\sum_{i=1}^{k}x_i<
\sum_{i=k+1}^{k+l}x_i$, then
 $\rho(G)<
\rho(G^*)$. 
\end{theorem}

\pf 
Then we have
\begin{align*}
x^T(A(G^*)-A(G))x
&=-x_v\sum_{i=1}^kx_i+x_v(\sum_{j=k+1}^{k+l}x_j-\sum_{i=1}^kx_i)+x_v\sum_{j=k+1}^{k+l}x_j\\
&=2x_v(\sum_{j=k+1}^{k+l}x_j-\sum_{i=1}^kx_i)\geq 0.
\end{align*}
So we have
\begin{align}\label{eq:11}
\rho(G^*) = \max_{\parallel y\parallel=1} y^TA(G^*)y\geq x^TA(G^*)x
\geq x^TA(G)x =\rho(G).
\end{align}
If $\rho(G^*) = \rho(G)$, then the equalities in (\ref{eq:11}) hold.
Thus
$$ x^TA(G^*)x=x^TA(G)x.$$
Hence $\sum_{i=k}^{k}x_i= \sum_{i=k+1}^{k+l}x_i$.  This completes
the proof. \qed

With above proof, we obtain the following result.

\begin{coro}
 Suppose $G^*$ in Theorem \ref{thm1} is connected, and $y = (y_1, y_2,
\ldots, y_n)^T$ is the Perron vector of $G^* $, then
$\sum_{i=k+1}^{k+l}y_i\geq \sum_{i=1}^{k}y_i$. Furthermore, if
$\sum_{i=k+1}^{k+l}x_i> \sum_{i=1}^{k}x_i$, then
$\sum_{i=k+1}^{k+l}y_i> \sum_{i=1}^{k}y_i$.
\end{coro}

\pf Suppose that $\sum_{i=k+1}^{k+l}y_i<\sum_{i=1}^{k}y_i$, by
Theorem \ref{thm1},  we have $\rho(G^*)<\rho(G)$, a contradiction.

Since $\sum_{i=k+1}^{k+l}x_i> \sum_{i=1}^{k}x_i$, by Theorem
\ref{thm1},  we have $\rho(G^*)>\rho(G)$. If
$\sum_{i=k+1}^{k+l}y_i\leq \sum_{i=1}^{k}y_i$, by Theorem \ref{thm1}
then we have $\rho(G^*)\leq\rho(G)$, a contradiction. This completes
the proof. \qed

\begin{lemma}\label{lem1}
If $\delta(G)> \frac{n+k}{2}+1$, then $G$ is $(k+1)$-connected.
\end{lemma}

\begin{theorem}
Let $n$, $k$ and $ \delta$ be three positive integers.  Among all
the connected graphs of order $n$ with connectivity at most $k$ and
minimum degree $ \delta$, the maximal spectral radius is obtained
uniquely at $G_{k,\delta,n}$.

\end{theorem}

\pf By Lemma \ref{lem1}, we have $2\delta\leq n+k+2$. If $n=k+1$,
then $K_{k+1}$ is an unique $k$-connected graph with order $n$. So
we can assume that $n\geq k+2$. Now we have to prove that for every
$G\in \mathcal {V}_{k,\delta,n} $, then $\rho(G)\leq
\rho(G_{k,\delta,n} )$, where the equality holds if and only if $G=
G_{k,\delta,n}$. Let $G^*\in \mathcal {V}_{k,\delta,n}$ with
$V(G^*)=\{v_1,\ldots,v_n\}$ be the graph with maximum spectral
radius in $\mathcal {V}_{k,\delta,n}$, that is, $\rho(G)\leq
\rho(G^*)$ for all $G\in \mathcal {V}_{k,\delta,n}$.

Denote the Perron vector with $x=(x_1,\ldots,x_n)$, where $x_i$
corresponding to $v_i$ for $i=1,\ldots,n$. Since $G^*\in \mathcal
{V}_{k,\delta,n}$ and it is not a complete graph, then $G^*$ has a
$k$-vertex cut, say $S=\{v_1,\ldots,v_k\}$. In the following, we
will prove the following three claims.

\medskip
\textbf{ Claim 1.~} $G^*$ contains exactly two components.
\medskip

Suppose contrary that $G^*-S$ contains three components $G_1$, $G_2$
and $G_3$. Let $u\in G_1$ and $v\in G_2$. It is obvious that $S$ is
also an $k$-vertex cut of $G+ uv$; i.e. $G^*+uv\in \mathcal
{V}_{k,\delta,n}$. By Theorem \ref{lem2}, we have $\rho(G^*) <
\rho(G^*+uv)$. This contradicts the definition of $G^*$.

Therefore, $G^*-S$ has exactly two components $G_1$ and $G_2$.

\medskip
\textbf{ Claim 2.~} Each subgraph of $G^*$ induced by vertices
$V(G_i)\cup S,$ for $  i = 1, 2$, is a clique.
\medskip

Suppose contrary that there is a pair of non-adjacent vertices $u,
v\in V(G_i)\cup S$ for $i = 1$ or $2$. Again, $G^*+ uv\in \mathcal
{V}_{k,\delta,n}$. By Theorem \ref{lem2}, we have $\rho(G^*) <(G^*+
uv)$. This contradicts the definition of $G^*$.

From Claim 2, it is clear that all $G_1$ and $G_2$ are cliques too.
Then we write $K_{n_i}$ instead of $G_i$, for $i = 1, 2$, in the
rest of the proof, where $n_i = |G_i|$. Since $\delta(G)\geq k $, we
have $n_i\geq   \delta-k $ for $i=1,2$.

\medskip
\textbf{ Claim 3.~} Either $n_1 = \delta-k+1$ or $n_2 = \delta-k+1$.
\medskip

Otherwise, we have $n_1 > \delta-k+1$ and $n_2 > \delta-k+1$. Let $v
\in G_1$ and $u\in G_2$ . Suppose    $$N_{G^*}(v) =
\{v_1,v_2,\ldots,v_{n_2-1}, v_1, v_2, \ldots, v_k \}$$ and
$$N_{G^*}(u) = \{u_1, u_2,
\ldots,u_{n_1-1},v_1, v_2, \ldots , v_k \}.$$

Partition the vertex set of $G$ into three parts: the   vertices of
$S$; the vertices of $G_1$; the   vertices of $G_2$. This is an
equitable partition of $G$ with quotient matrix
\[Q=\left(\begin{matrix}
 k-1  & n_1 & n_2 \\
k     & n_1 -1 &
0\\
k     & 0 & n_2-1
\end{matrix} \right)
\]
By Perron-Frobenius Theorem, $Q$ has a Perron-vector
$x=\{x_1,x_2,x_3\}$. Now we show that  $x_2<x_3$ if $n_1<n_2$. Let
$\rho(Q)$ denotes the largest eigenvalue of $Q$. Then we have
\begin{align}
k x_1+   ( n_1-1)x_2 =\rho(Q)x_2\label{eq21}\\
k x_1+   ( n_2-1)x_3 =\rho(Q)x_3\label{eq22}
\end{align}
By (\ref{eq21}) and (\ref{eq22}), we have
\begin{align*}
( n_2-1)x_3-( n_1-1)x_2 =\rho(Q)(x_3-x_2).
\end{align*}
Hence
\begin{align*}
( \rho(Q)-n_2+1)(x_3-x_2)=( n_2-n_1)x_2 >0.
\end{align*}
Since $\rho(Q)$ is also the largest eigenvalue of $G^*$,   we have
$\rho(Q)> \rho(K_{n_2})=n_2-1$. Hence $x_3>x_2$. The eigenvector $x$
can be extended to an  eigenvector of $A(G^*)$, say
$$y=(x_{11},\ldots,x_{1k},x_{21},\ldots,x_{2n_1},x_{31},\ldots,x_{3n_2}),$$
 where $x_{i1}=\ldots=x_{in_i}=x_i$ for $i=2,3$ and
$x_{11}=\cdots=x_{1k}=x_1$. Let
$z=\frac{1}{\sqrt{kx_1^2+n_1x_2^2+n_3x_3^2}}y$. We have $zz^T=1$ and
so $z$ is a Perron-vector of $G^*$. Let $G = G^*- \{vv_1,vv_2,
\ldots, vv_k\} + \{vv_{k+1},\ldots, vv_{k+l}\}$ and we have $G\in
\mathcal {V}_{k,\delta,n}$. Since $n_2x_3>(n_1-1)x_2$, by Theorem
\ref{thm1}, $\rho(G^*) < \rho(G)$, which is a contradiction. This
completes claim 3.

By Claim 3, we have  $n_1=\delta+1-k$. Hence $G^*=G_{k,\delta,n}$.
This completes the proof. \qed

\begin{theorem}
The spectral radius of $G_{k,\delta,n}$ is the largest root of the
following equation
$x^3+(3-n)x^2+(n\delta-\delta^2-n-kn+k+k\delta+2-2\delta)x+(kn\delta+k^2+n\delta+k^2\delta-
k\delta-k^2n-k\delta^2-2\delta-\delta^2)=0.$
\end{theorem}

\pf Let $G_{1}$ be the subgraph of $G_{k,\delta,n}$ induced by $k$
vertices of all the vertices of degree $n-1$, $G_{2}$ be the
subgraph induced by all the vertices of degree $\delta$ vertices and
$G_{3}$ be the subgraph induced by the remaining $n-\delta-1$
vertices. Also, let $G_{ij}$ be the bipartite subgraph induced by
$V(G_i)$ and $V(G_j)$ and let $e_{ij}$ be the size of $G_{12}$. A
theorem of Haemers \cite{Haemers} shows that eigenvalues of the
quotient matrix of the partition interlace the eigenvalues of the
adjacency matrix of $G$. The quotient matrix $Q$ is the following
\[Q=\left(\begin{matrix}
\frac{2e_{1}}{n_{1}}   & \frac{e_{12}}{n_{1}} & \frac{e_{13}}{n_{1}} \\
\frac{e_{21}}{n_{2}}     & \frac{2e_{2}}{n_{2}}&
\frac{e_{23}}{n_{2}}\\
\frac{e_{31}}{n_{3}}     & \frac{e_{32}}{n_{3}}&
\frac{2e_{33}}{n_{3}}
\end{matrix} \right)\\
=\begin{pmatrix}
k-1   & \delta-k+1 & n-\delta-1  \\
k    & \delta-k & 0 \\
k & 0 & n-\delta-2
\end{pmatrix}.
\]
Applying eigenvalue interlacing  to the greatest eigenvalue of $G$,
we get
\begin{align}\label{eq:4}
\lambda_{1}(H)\geq \lambda_{1}(Q),
\end{align}
with the equality  if the partition is equitable [\cite{Godsil},
p.202]. Note that the partition is equitable, so the equality hold.
This completes the proof.  \qed

\end{document}